# CENTRAL LIMIT THEOREM FOR SEQUENTIAL MONTE CARLO METHODS AND ITS APPLICATION TO BAYESIAN INFERENCE

### By Nicolas Chopin

### *Bristol University*

The term "sequential Monte Carlo methods" or, equivalently, "particle filters," refers to a general class of iterative algorithms that performs Monte Carlo approximations of a given sequence of distributions of interest $(\pi_t)$. We establish in this paper a central limit theorem for the Monte Carlo estimates produced by these computational methods. This result holds under minimal assumptions on the distributions $\pi_t$, and applies in a general framework which encompasses most of the sequential Monte Carlo methods that have been considered in the literature, including the resample-move algorithm of Gilks and Berzuini [*J. R. Stat. Soc. Ser. B Stat. Methodol.* **63** (2001) 127–146] and the residual resampling scheme. The corresponding asymptotic variances provide a convenient measurement of the precision of a given particle filter. We study, in particular, in some typical examples of Bayesian applications, whether and at which rate these asymptotic variances diverge in time, in order to assess the long term reliability of the considered algorithm.

**1. Introduction.** Sequential Monte Carlo methods form an emerging, yet already very active branch of the Monte Carlo paradigm. Their growing popularity comes in part from the fact that they are often the only viable computing techniques in those situations where data must be processed sequentially. Their range of applicability is consequently very wide, and includes nonexclusively signal processing, financial modeling, speech recognition, computer vision, neural networks, molecular biology and genetics, target tracking and geophysics, among others. A very good introduction to the field has been written by Künsch (2001), while the edited volume of Doucet, de Freitas and Gordon (2001) provides an interesting coverage of recent developments in theory and applications.

Specifically, sequential Monte Carlo methods (alternatively termed "particle filters" or "recursive Monte Carlo filters") are iterative algorithms that produce and update recursively a set of weighted simulations (the "particles") in order to provide a Monte Carlo approximation of an evolving distribution of interest $\pi_t(d\theta_t)$, $t$ being an integer index. In a sequential Bayesian framework, $\pi_t(d\theta_t)$ will usually represent the posterior distribution of parameter $\theta_t$ given the $t$ first observations. The term "parameter" must be understood here in a broad sense, in that $\theta_t$ may include any unknown quantity which may be inferred from the $t$ first observations, and is not necessarily of constant dimension. We denote by $\Theta_t$ the support of $\pi_t(d\theta_t)$.

The study of the asymptotic properties of sequential Monte Carlo methods is admittedly a difficult problem, and some methodological papers [Liu and Chen (1998), e.g.] simply state some form of the law of large numbers for the most elaborate algorithms, that is, the Monte Carlo estimates are shown to converge almost surely to the quantity of interest as $H$, the number of particles, tends toward infinity. More refined convergence results have been obtained, such as the central limit theorem of Del Moral and Guionnet (1999), later completed by Del Moral and Miclo (2000), or upper bounds for the Monte Carlo error expressed in various norms [Crisan and Lyons (1997, 1999), Crisan, Gaines and Lyons (1998), Crisan and Doucet (2000), Del Moral and Guionnet (2001), Künsch (2001) and Le Gland and Oudjane (2004)]. Unfortunately, it has been, in general, at the expense of generality [with the exception of Crisan and Doucet (2000)], whether in terms of computational implementation (only basic algorithms are considered, which may not be optimal) or of applicability (the sequence $\pi_t$ has to be generated from some specific dynamical model that fulfills various conditions).

In this paper we derive a central limit theorem that applies to most of the sequential Monte Carlo techniques developed recently in the methodological literature, including the resample-move algorithm of Gilks and Berzuini (2001), the auxiliary particle filter of Pitt and Shephard (1999) and the stochastic remainder resampling scheme [Baker (1985, 1987)], also known as the residual resampling scheme [Liu and Chen (1998)]. No assumption is made on the model that generates the sequence of distributions of interest $(\pi_t)$, so that our theorem equally applies to those recent algorithms [Chopin (2002), Del Moral and Doucet (2002) and Cappé, Guillin, Marin and Robert (2004)] that have been developed for contexts that widely differ from the standard application of sequential Monte Carlo methods, namely, the sequential analysis of state space models.

The appeal of a central limit theorem is that it provides an (asymptotically) exact measure of the Monte Carlo error, through the asymptotic variance. This allows for a rigorous comparison of the relative efficiency of given algorithms. In this way, we show in this paper, again by comparing the appropriate asymptotic variances, that the residual resampling scheme always outperforms the multinomial resampling scheme, and that the Rao–Blackwell variance reduction technique of Doucet, Godsill and Andrieu (2000) is, indeed, effective.

The most promising application of our central limit theorem is the possibility to assess the stability of a given particle filter (in terms of precision of the computed estimates) through the time behavior of the corresponding asymptotic variances. This is a critical issue since it is well known that sequential Monte Carlo methods tend to degenerate in a number of cases, sometimes at a very fast rate. We consider in this paper some typical Bayesian problems, such as the sequential analysis of state-space models. We will show that under some conditions stability can be achieved at least for "filtering" the states, that is, for approximating the marginal posterior density $\pi_t(x_t)$, where $x_t$ stands for the current state at iteration $t$.

The paper is organized as follows. Section 2 proposes a generic description of particle filters, establishes a central limit theorem for computed estimates in a general framework and draws some conclusions from this result. Section 3 discusses the stability of particle filters through the time behavior of the asymptotic variances provided by the central limit theorem. Proofs of theorems are put in the Appendix.

## 2. Central limit theorem for particle filters.

2.1. *General formulation of particle filters.* In full generality, a particle system is a triangular array of random variables in $\Theta \times \mathbb{R}^+$,

$$(\theta^{(j,H)}, w^{(j,H)})_{j \leq H},$$

where $\Theta$ is some space of interest. The variables $\theta^{(j,H)}$ are usually called "particles," and their contribution to the sample may vary according to their weights $w^{(j,H)}$. We will say that this particle system *targets* a given distribution $\pi$ defined on $\Theta$ if and only if

$$(1) \qquad \frac{\sum_{j=1}^H w^{(j,H)} \varphi(\theta^{(j,H)})}{\sum_{j=1}^H w^{(j,H)}} \to \mathbb{E}_\pi(\varphi)$$

holds almost surely as $H \to +\infty$ for any measurable function $\varphi$ such that the expectation above exists. A first example of a particle system is a denumerable set of independent draws from $\pi$, with unit weights, which obviously targets $\pi$. In this simple case, particles and weights do not depend on $H$, and the particle system is a sequence rather than a triangular array. This is not the case in general, however, and, while cumbersome, the dependence in $H$ will be maintained in notation to allow for a rigorous mathematical treatment.

Now assume a sequence $(\pi_t)_{t \in \mathbb{N}}$ of distributions defined on a sequence of probabilized spaces $(\Theta_t)$. In most, if not all, applications, $\Theta_t$ will be a power of the real line or some subset of it, and, henceforth, $\pi_t(\cdot)$ will also denote the density of $\pi_t$ with respect to an appropriate version of the Lebesgue measure. A sequential Monte Carlo algorithm (or particle filter) is a method for producing a particle system whose target evolves in time: at iteration $t$ of the algorithm, the particle system targets $\pi_t$, and therefore allows for Monte Carlo approximations of the distribution of (current) interest $\pi_t$. Clearly, particle filters do not operate in practice on infinite triangular arrays but rather manipulate particle vectors of fixed size $H$. One must keep in mind, however, that the justification of such methods is essentially asymptotic.

The structure of a particle filter can be decomposed into three basic iterative operations, that will be referred to hereafter as mutation, correction and selection steps. At the beginning of iteration $t$, consider a particle system $(\hat{\theta}_{t-1}^{(j,H)}, 1)_{j \leq H}$, that is, with unit weights, which targets $\pi_{t-1}$. The mutation step consists in producing new particles drawn from

$$\theta_t^{(j,H)} \sim k_t(\hat{\theta}_{t-1}^{(j,H)}, d\theta_t),$$

where $k_t$ is a transition kernel which maps $\Theta_{t-1}$ into $\mathcal{P}(\Theta_t)$, the set of probability measures on $\Theta_t$. The "mutated" particles (with unit weights) target the new distribution $\tilde{\pi}_t(\cdot) = \int \pi_{t-1}(\theta_{t-1}) k_t(\theta_{t-1}, \cdot) \, d\theta_{t-1}$. This distribution $\tilde{\pi}_t$ is usually not relevant to the considered application, but rather serves as an intermediary stage for practical reasons. To shift the target to the distribution of interest $\pi_t$, particles are assigned weights

$$w_t^{(j,H)} \propto v_t(\theta_t^{(j,H)}) \qquad \text{with } v_t(\cdot) = \pi_t(\cdot)/\tilde{\pi}_t(\cdot).$$

This is the correction step. The particle system $(\theta_t^{(j,H)}, w_t^{(j,H)})_{j \leq H}$ targets $\pi_t$. The function $v_t$ is referred to as the weight function. Note that the normalizing constants of the densities $\pi_t$ and $\tilde{\pi}_t$ are intractable in most applications. This is why weights are defined up to a multiplicative constant, which has no bearing anyway on the estimates produced by the algorithm, since they are weighted averages.

Finally, the selection step consists in replacing the current vector of particles by a new, uniformly weighted vector $(\hat{\theta}_t^{(j,H)}, 1)_{j \leq H}$, which contains a number $n^{(j,H)}$ of replicates of particle $\theta_t^{(j,H)}$, $n^{(j,H)} \geq 0$. The $n^{(j,H)}$'s are random variables such that $\sum_j n^{(j,H)} = H$ and $\mathbb{E}(n^{(j,H)}) = H\rho_j$, where the normalized weights are given by

$$\rho_j = w_t^{(j,H)} \Big/ \sum_{j=1}^{H} w_t^{(j,H)},$$

and where dependencies in $H$ and $t$ are omitted for convenience. In this way, particles whose weights are too small are discarded, while particles with important weights serve as multiple starting points for the next mutation step. There are various ways of generating the $n^{(j,H)}$'s. Multinomial resampling [Gordon, Salmond and Smith (1993)] amounts to drawing independently the $H$ new particles from the multinomial distribution which produces $\theta_t^{(j,H)}$ with probability $\rho_j$. Residual resampling [originally termed "stochastic remainder sampling" in the genetic algorithm literature, Baker (1985, 1987), then rediscovered by Liu and Chen (1998)] consists in reproducing $\lfloor H\rho_j \rfloor$ times each particle $\theta_t^{(j,H)}$, where $\lfloor \cdot \rfloor$ stands for the integer part. The particle vector is completed by $H^r = H - \sum_j \lfloor H\rho_j \rfloor$ independent draws from the multinomial distribution which produces $\theta_t^{(j,H)}$ with probability $(H\rho_j - \lfloor H\rho_j \rfloor)/H^r$. Systematic resampling [another method initially proposed in the genetic algorithm field, Whitley (1994), then rediscovered by Carpenter, Clifford and Fearnhead (1999); see also Crisan and Lyons (2002) for a slightly different algorithm] is another interesting selection scheme, which is such that the number of replicates $n^{(j,H)}$ is ensured to differ from $H\rho_j$ by at most one. We failed, however, to extend our results to this third selection scheme.

The structure of a particle filter can be summarized as follows:

1. Mutation: Draw for $j = 1, \ldots, H$,

$$\theta_t^{(j,H)} \sim k_t(\hat{\theta}_{t-1}^{(j,H)}, d\theta_t),$$

   where $k_t : \Theta_{t-1} \to \mathcal{P}(\Theta_t)$ is a given probability kernel.

2. Correction: Assign weights to particles so that, for $j = 1, \ldots, H$,

$$w_t^{(j,H)} \propto v_t(\theta_t^{(j,H)}) = \pi_t(\theta_t^{(j,H)})/\tilde{\pi}_t(\theta_t^{(j,H)}),$$

   where $\tilde{\pi}_t(\cdot) = \int \pi_{t-1}(\theta_{t-1}) k_t(\theta_{t-1}, \cdot) \, d\theta_{t-1}$.

3. Selection: Resample, according to a given selection scheme,

$$(\theta_t^{(j,H)}, w_t^{(j,H)})_{j \le H} \to (\hat{\theta}_t^{(j,H)}, 1)_{j \le H}.$$

The first mutation step, $t = 0$, is assumed to draw independent and identically distributed particles from some instrumental distribution $\tilde{\pi}_0$.

It is shown without difficulty that the particle system produced by this generic algorithm does iteratively target the distributions of interest, that is, the following convergences hold almost surely:

$$H^{-1} \sum_{j=1}^{H} \varphi(\theta_t^{(j,H)}) \to \mathbb{E}_{\tilde{\pi}_t}(\varphi),$$

$$\frac{\sum_{j=1}^{H} w_t^{(j,H)} \varphi(\theta_t^{(j,H)})}{\sum_{j=1}^{H} w_t^{(j,H)}} \to \mathbb{E}_{\pi_t}(\varphi),$$

$$H^{-1} \sum_{j=1}^{H} \varphi(\hat{\theta}_t^{(j,H)}) \to \mathbb{E}_{\pi_t}(\varphi),$$

as $H \to +\infty$, provided these expectations exist. These convergences will be referred to as the law of large numbers for particle filters.

2.2. *Some examples of particle filters.* The general formulation given in the previous section encompasses most of the sequential Monte Carlo algorithms described in the literature. By way of illustration, assume first that the distributions $\pi_t$ are defined on a common space, $\Theta_t = \Theta$. In a Bayesian framework, $\pi_t$ will usually be the posterior density of $\theta$, given the $t$ first observations, $\pi_t(\theta) = \pi(\theta|y_{1:t})$, where $y_{1:t}$ denotes the sequence of observations $y_1, \ldots, y_t$. If particles are not mutated, $k_t$ being the "identity kernel" $k_t(\theta, \cdot) = \delta_\theta$, we have $\tilde{\pi}_t = \pi_{t-1}$ for $t > 0$, and our generic particle filter becomes one of the variations of the sequential importance resampling algorithm [Rubin (1988), Gordon, Salmond and Smith (1993) and Liu and Chen (1998)]. The weight function simplifies to

$$\upsilon_t(\theta) = \pi(\theta|y_{1:t})/\pi(\theta|y_{1:t-1}) \propto p(y_t|y_{1:t-1}, \theta)$$

in a Bayesian model, where $p(y_t|y_{1:t-1}, \theta)$ is the conditional likelihood of $y_t$, given the parameter $\theta$ and previous observations.

Gilks and Berzuini (2001) propose a variant of this algorithm, namely, the resample-move algorithm, in which particles are mutated according to an MCMC [Markov chain Monte Carlo; see, e.g., Robert and Casella (1999)] kernel $k_t$, which admits $\pi_{t-1}$ as an invariant density. In that case, we still have $\tilde{\pi}_t = \pi_{t-1}$, and the expression for the weight function $\upsilon_t$ is unchanged. The motivation of this strategy is to add new particle values along iterations so as to limit the depletion of the particle system.

Now consider the case where $\pi_t$ is defined on a space of increasing dimension of the form $\Theta_t = \mathcal{X}^t$. A typical application is the sequential inference of a dynamical model which involves a latent process $(x_t)$, and $\pi_t$ stands then for density $\pi(x_{1:t}|y_{1:t})$. Assume $k_t$ can be decomposed as

$$k_t(x_{1:t-1}^*, dx_{1:t}) = \kappa_t(x_{1:t-1}^*, dx_{1:t-1})q_t(x_t|x_{1:t-1})\,dx_t,$$

where $\kappa_t : \mathcal{X}^{t-1} \to \mathcal{P}(\mathcal{X}^{t-1})$ is a transition kernel, and $q_t(\cdot|\cdot)$ is some conditional probability density. If $\kappa_t$ admits $\pi_{t-1}$ as an invariant density, the weight function is given by

$$(2) \qquad \upsilon_t(x_{1:t}) = \frac{\pi_t(x_{1:t})}{\pi_{t-1}(x_{1:t-1})q_t(x_t|x_{1:t-1})}.$$

Again, the case where $\kappa_t$ is the identity kernel corresponds to some version of the sequential importance resampling algorithm, while setting $\kappa_t$ to a given

MCMC transition kernel with invariant density $\pi_{t-1}$ leads to the resample-move algorithm of Gilks and Berzuini (2001). The standard choice for $q_t(\cdot|\cdot)$ is the conditional prior density of $x_t$, given $x_{1:t-1}$, as suggested originally by Gordon, Salmond and Smith (1993), but this is not always optimal, as pointed out by Pitt and Shephard (1999) and Doucet, Godsill and Andrieu (2000). In fact, it is generally more efficient to build some conditional density $q_t$ which takes into account the information carried by $y_t$ in some way, in order to simulate more values compatible with the observations.

These two previous cases can be combined into one, by considering a dynamical model which features at the same time a fixed parameter $\theta$ and a sequence of latent variables $(x_t)$, so that $\Theta_t = \Theta \times \mathcal{X}^t$, and $\pi_t$ stands for the joint posterior density $\pi(\theta, x_{1:t}|y_{1:t})$.

2.3. *Central limit theorem.* The following quantities will play the role of asymptotic variances in our central limit theorem. Let, for any measurable $\varphi : \Theta_0 \to \mathbb{R}^d$, $\widetilde{V}_0(\varphi) = \operatorname{Var}_{\tilde{\pi}_0}(\varphi)$, and by induction, for any measurable $\varphi : \Theta_t \to \mathbb{R}^d$,

$$\widetilde{V}_t(\varphi) = \widehat{V}_{t-1}\{\mathbb{E}_{k_t}(\varphi)\} + \mathbb{E}_{\pi_{t-1}}\{\operatorname{Var}_{k_t}(\varphi)\}, \qquad t > 0, \tag{3}$$

$$V_t(\varphi) = \widetilde{V}_t\{v_t \cdot (\varphi - \mathbb{E}_{\pi_t}\varphi)\}, \qquad t \geq 0, \tag{4}$$

$$\widehat{V}_t(\varphi) = V_t(\varphi) + \operatorname{Var}_{\pi_t}(\varphi), \qquad t \geq 0. \tag{5}$$

The notation $\mathbb{E}_{k_t}(\varphi)$ and $\operatorname{Var}_{k_t}(\varphi)$ is shorthand for the functions $\mu(\theta_{t-1}) = \mathbb{E}_{k_t(\theta_{t-1}, \cdot)}\{\varphi(\cdot)\}$ and $\Sigma(\theta_{t-1}) = \operatorname{Var}_{k_t(\theta_{t-1}, \cdot)}\{\varphi(\cdot)\}$, respectively. Note that these equations do not necessarily produce finite variances for any $\varphi$. We now specify the classes of functions for which the central limit theorem enunciated below will hold, and, in particular, for which these asymptotic variances exist. Denoting by $\|\cdot\|$ the Euclidean norm in $\mathbb{R}^d$, we define recursively $\Phi_t^{(d)}$ to be the set of measurable functions $\varphi : \Theta_t \to \mathbb{R}^d$ such that for some $\delta > 0$,

$$\mathbb{E}_{\tilde{\pi}_t}\|v_t \cdot \varphi\|^{2+\delta} < +\infty, \tag{6}$$

and that the function $\theta_{t-1} \mapsto \mathbb{E}_{k_t(\theta_{t-1}, \cdot)}\{v_t(\cdot)\varphi(\cdot)\}$ is in $\Phi_{t-1}^{(d)}$. The initial set $\Phi_0^{(d)}$ contains all the measurable functions whose moments of order two with respect to $\tilde{\pi}_0$ are finite.

THEOREM 1. *If the selection step consists of multinomial resampling, and provided that the unit function $\theta_t \mapsto 1$ belongs to $\Phi_t^{(1)}$ for every $t$, then for any $\varphi \in \Phi_t^{(d)}$, $\mathbb{E}_{\pi_t}(\varphi)$, $V_t(\varphi)$ and $\widehat{V}_t(\varphi)$ are finite quantities, and the following convergences in distribution hold as $H \to +\infty$:*

$$H^{1/2}\left\{\frac{\sum_{j=1}^H w_t^{(j,H)}\varphi(\theta_t^{(j,H)})}{\sum_{j=1}^H w_t^{(j,H)}} - \mathbb{E}_{\pi_t}(\varphi)\right\} \xrightarrow{\mathcal{D}} \mathcal{N}\{0, V_t(\varphi)\},$$

$$H^{1/2}\left\{H^{-1}\sum_{j=1}^{H}\varphi(\hat{\theta}_t^{(j,H)}) - \mathbb{E}_{\pi_t}(\varphi)\right\} \xrightarrow{\mathcal{D}} \mathcal{N}\{0, \widehat{V}_t(\varphi)\}.$$

A proof is given in the Appendix. In the course of the proof an additional central limit theorem is established for the unweighted particle system $(\theta_t^{(j,H)}, 1)$ produced by the mutation step, which targets $\tilde{\pi}_t$. This result is not given here, however, for it holds for a slightly different class of functions, and is of less practical interest. The assumption that the function $\theta_t \mapsto 1$ belongs to $\Phi_t^{(1)}$ deserves further comment. Qualitatively, it implies that the weight function $v_t$ has finite moment of order $2 + \delta$ with respect to $\tilde{\pi}_t$, for some $\delta > 0$, and, therefore, restricts somehow the dispersion of the particle weights. It also implies that $\Phi_t^{(d)}$ contains all bounded functions $\varphi$. In practice this assumption will be fulfilled, for instance, whenever each weight function $v_t$ is bounded from above, which occurs in many practical settings.

A central limit theorem also holds when the selection step follows the residual sampling scheme of Liu and Chen (1998), but this imposes some change in the expression for the asymptotic variances. The new expression for $\widehat{V}_t(\varphi)$ is

(7) $$\widehat{V}_t(\varphi) = V_t(\varphi) + R_t(\varphi),$$

where

(8) $$R_t(\varphi) = \mathbb{E}_{\tilde{\pi}_t}\{r(v_t)\varphi\varphi'\} - \frac{1}{\mathbb{E}_{\tilde{\pi}_t}\{r(v_t)\}}[\mathbb{E}_{\tilde{\pi}_t}\{r(v_t)\varphi\}][\mathbb{E}_{\tilde{\pi}_t}\{r(v_t)\varphi\}]',$$

and $r(x)$ is $x$ minus its integer part.

THEOREM 2. *The results of Theorem 1 still hold when the selection steps consists of residual resampling, except that the asymptotic variances are now defined by equations* (3), (4) *and* (7).

The proofs of Theorems 1 and 2 (given in the Appendix) rely on an induction argument: conditional on past iterations, each step generates independent (but not identically distributed) particles, which follow some (conditional) central limit theorem. In contrast, the systematic resampling scheme is such that, given the previous particles, the new particle system is entirely determined by a single draw from a uniform distribution; see Whitley (1994). This is why extending our results to this third selection scheme seems not straightforward, and possibly requires an entirely different approach.

The appeal of the recursive formulae (3)–(5) and (7) is that they put forward the impact of each new step on the asymptotic variance, particularly the additive effect of the selection and mutation steps. In the multinomial case, an alternative expression for the asymptotic variance is

$$V_t(\varphi) = \sum_{k=0}^{t} \mathbb{E}_{\tilde{\pi}_k}[v_k^2 \mathcal{E}_{k+1:t}\{\varphi - \mathbb{E}_{\pi_t}(\varphi)\} \mathcal{E}_{k+1:t}\{\varphi - \mathbb{E}_{\pi_t}(\varphi)\}'], \tag{9}$$

where $\mathcal{E}_t$ is the functional operator which associates to $\varphi$ the function

$$\mathcal{E}_t(\varphi) : \theta_{t-1} \mapsto \mathbb{E}_{k_t(\theta_{t-1}, \cdot)}\{v_t(\cdot)\varphi(\cdot)\}, \tag{10}$$

and $\mathcal{E}_{k+1:t}(\varphi) = \mathcal{E}_{k+1} \circ \cdots \circ \mathcal{E}_t(\varphi)$ for $k+1 \leq t$, $\mathcal{E}_{t+1:t}(\varphi) = \varphi$. This closed form expression is more convenient when studying the stability of the asymptotic variance over time, as we will illustrate in the next section. A similar formula for the residual case can be obtained indirectly by deriving the difference between the multinomial and the residual cases, that is, for $t > 0$,

$$V_t^r(\varphi) - V_t(\varphi) = \sum_{k=0}^{t-1} [R_k\{\mathcal{E}_{k+1:t}(\varphi)\} - \mathrm{Var}_{\pi_k}\{\mathcal{E}_{k+1:t}(\varphi)\}], \tag{11}$$

where $V_t(\varphi)$, $V_t^r(\varphi)$ are defined through the recursions (3)–(5) and (3), (4) and (7), respectively. In the following, we will similarly distinguish the residual case through an $r$-suffix in notation.

2.4. *First conclusions.* A first application of this central limit theorem is to provide a rigorous justification for some heuristic principles that have been stated in the literature, see, for instance, Liu and Chen (1998). Inequalities in this section refer to the canonical order for symmetric matrices, that is to say $A > B$ (resp. $A \geq B$) if and only if $A - B$ is positive definite (resp. positive semidefinite).

First, it is preferable to compute any estimate before the selection step, since the immediate effect of the latter is a net increase in asymptotic variance: $\widehat{V}_t(\varphi) > V_t(\varphi)$ for any nonconstant function $\varphi$. In this respect one may wonder why selection steps should be performed. We will see that the immediate degradation of the particle system is often largely compensated for by gains in precision in the future iterations.

Second, residual sampling always outperforms multinomial resampling. Let $\varphi : \Theta_t \to \mathbb{R}^d$ and $\bar{\varphi} = \varphi - \mathbb{E}_{\pi_t}(\varphi)$. Then

$$R_t(\varphi) = R_t(\bar{\varphi}) \leq \mathbb{E}_{\tilde{\pi}_t}\{r(v_t)\bar{\varphi}\bar{\varphi}'\} \leq \mathrm{Var}_{\pi_t}(\varphi),$$

since $r(x) \leq x$. It follows from this inequality and (11) that $V_t^r(\varphi) \leq V_t(\varphi)$. Actually, a substantial gain should be expected when using the residual scheme since the inequality above is clearly not sharp.

Our central limit theorem also provides a formal justification for resorting to "marginalized" particle filters, as explained in the following section.

2.5. *Marginalized particle filters.* In some specific cases it is possible to decompose the density $\pi_t(\theta_t)$ into $\pi_t^m(\xi_t)\pi_t^c(\lambda_t|\xi_t)$, with $\theta_t = (\xi_t, \lambda_t)$ lying in $\Theta_t = \Xi_t \times \Lambda_t$, in such a way that it is possible to implement a particle filter that targets the marginal densities $\pi_t^m$ rather than the $\pi_t$'s. When this occurs, this second algorithm usually produces more precise estimators (in a sense that we explain below) in the $\xi_t$-dimension. The idea of resorting to "marginalized" particle filters has been formalized by Doucet, Godsill and Andrieu (2000), and implemented in various settings by Chen and Liu (2000), Chopin (2001) and Andrieu and Doucet (2002), among others.

Doucet, Godsill and Andrieu's (2000) justification for resorting to "marginalized" particle filters is that they yield importance weights with a smaller variance than their "unmarginalized" counterpart, which suggests that the produced estimates are also less variable. This is proven by a Rao–Blackwell decomposition, and, consequently, "marginalized" particle filters are sometimes referred to as "Rao–Blackwellized" particle filters. We now extend the argument of these authors by proving that the asymptotic variance of any estimator is, indeed, smaller in the "marginalized" case. Assume decompositions of $\pi_t$ and $\tilde{\pi}_t$ of the form

$$\pi_t(\theta_t) = \pi_t^m(\xi_t)\pi_t^c(\lambda_t|\xi_t), \qquad \tilde{\pi}_t(\theta_t) = \tilde{\pi}_t^m(\xi_t)\tilde{\pi}_t^c(\lambda_t|\xi_t),$$

where $(\xi_t, \lambda_t)$ identifies to $\theta_t$, and $\pi_t^m$, $\pi_t^c$, $\tilde{\pi}_t^m$, $\tilde{\pi}_t^c$, are, respectively, marginal and conditional densities of $\xi_t$ and $\lambda_t$. Consider two particle filters, tracking, respectively, $(\pi_t)$ and $(\pi_t^m)$. It is assumed that both filters implement the same selection scheme (whether multinomial or residual), and that their mutation steps consist in drawing, respectively, from kernels $k_t$ and $k_t^m$, which are such that the following probability measures coincide on $\Theta_t = \Xi_t \times \Lambda_t$,

$$(12)\quad \begin{aligned}\int_{\Lambda_{t-1}} & \pi_{t-1}^c(\lambda_{t-1}|\xi_{t-1})k_t\{(\xi_{t-1},\lambda_{t-1}),(d\xi_t,d\lambda_t)\}\,d\lambda_{t-1} \\ & = k_t^m(\xi_{t-1},d\xi_t)\tilde{\pi}_t^c(\lambda_t|\xi_t)\,d\lambda_t,\end{aligned}$$

for almost every $\xi_{t-1}$ in $\Xi_{t-1}$. Note that in full generality it is not always possible to build a kernel $k_t^m$ from a given $k_t$ which satisfies this relation. As illustrated by the aforementioned references, however, it is feasible in some cases of interest. This equality implies, in particular, that

$$\int \pi_{t-1}^m(\xi_{t-1})k_t^m(\xi_{t-1},\cdot)\,d\xi_{t-1} = \tilde{\pi}_t^m(\cdot).$$

Asymptotic variances and other quantities are distinguished similarly through the $m$-suffix for the marginal case, that is, $V_t(\varphi)$ and $V_t^m(\varphi)$, and so on.

THEOREM 3.   *For any $\varphi : \Xi_t \to \mathbb{R}^d$ such that $\varphi \in \Phi_t^{(d)}$, we have $V_t^m(\varphi) \leq V_t(\varphi)$ and $V_t^{m,r}(\varphi) \leq V_t^r(\varphi)$. These inequalities are attained for a nonconstant $\varphi$ if and only if $\pi_t^c(\cdot|\xi_t) = \tilde{\pi}_t^c(\cdot|\xi_t)$ for almost every $\xi_t \in \Xi_t$, for any $t \geq 0$.*

As suggested by the condition for equality above or more clearly exhibited in the proof in the Appendix, marginalizing allows for canceling the weight dispersion due to the discrepancy between conditional densities $\tilde{\pi}_t^c$ and $\pi_t^c$, while the part due to the discrepancy between marginal densities $\pi_t^m$ and $\tilde{\pi}_t^m$ remains identical.

Beyond the small number of cases where this marginalization technique can be effectively carried out, this result has also strong qualitative implications. In the following sections we will study the behavior of the time sequence $V_t(\varphi)$ in order to measure whether and at which rate a given particle filter "diverges." In this respect, we will be able in some cases to build a marginalized particle filter whose rate of divergence is theoretically known, thus providing a lower bound for the actual rate of divergence of the considered particle filter.

## 3. Stability of particle filters.

3.1. *Sequential importance sampling.*   The sequential importance sampling algorithm is a particle filter that alternates mutation and correction steps, but does not perform any selection step. Weights are consequently not initialized to one at each iteration, and are rather updated through

$$w_t^{(j)} \propto w_{t-1}^{(j)} v_t(\theta_t^{(j)}).$$

We suppress any notational dependence on $H$ since it is meaningless in such a case. Due to its specific nature, this algorithm needs to be treated separately. Since particles are not resampled, they remain independent through iterations. It follows via the standard central limit theorem that

$$H^{1/2} \left\{ \frac{\sum_{j=1}^H w_t^{(j)} \varphi(\theta_t^{(j)})}{\sum_{j=1}^H w_t^{(j)}} - \mathbb{E}_{\pi_t}(\varphi) \right\} \xrightarrow{\mathcal{D}} \mathcal{N}\{0, V_t^{\mathrm{sis}}(\varphi)\},$$

where the corresponding asymptotic variance is

$$V_t^{\mathrm{sis}}(\varphi) = \mathbb{E}_{\tilde{\pi}_t} \left[ \frac{\pi_t}{\tilde{\pi}_t} \{\varphi - \mathbb{E}_{\pi_t}(\varphi)\} \right]^2,$$

and $\tilde{\pi}_t$ denotes this time the generating distribution of particles $\theta_t^{(j)}$ obtained by the recursion of mutation kernels $k_t(\cdot, \cdot)$, that is,

$$\tilde{\pi}_t(\cdot) = \int \tilde{\pi}_{t-1}(\theta_{t-1}) k_t(\theta_{t-1}, \cdot) \, d\theta_{t-1},$$

the distribution $\tilde{\pi}_0$ being arbitrary. Sequential importance sampling is rarely an efficient algorithm, but the value of $V_t^{\mathrm{sis}}(\varphi)$ can serve as a benchmark in some occasions, as we will see in the following.

3.2. *Sequential importance sampling and resampling in the fixed parameter case.* In the fixed parameter case, that is, $\Theta_t = \Theta$ and $\pi_t(\theta) = \pi(\theta|y_{1:t})$, $\pi_t$ is expected to become more and more informative on $\theta$, and to eventually converge to a Dirac mass at some point $\theta_0$. Sequential importance sampling and resampling algorithms typically diverge in such a situation, since they generate once and for all the set of particle values from $\tilde{\pi}_0$, a majority of which are presumably far from $\theta_0$. The following result quantifies this degeneracy effect.

THEOREM 4. *Let $\varphi : \Theta \to \mathbb{R}^d$, $\varphi \in \Phi_t^{(d)}$. Then under regularity conditions given in the Appendix, there exist positive constants $c_1$, $c_2$ and $c_3$ such that*

$$\|V_t^{\mathrm{sis}}(\varphi)\| \asymp c_1 t^{p/2-1}, \qquad \|V_t^r(\varphi)\| \asymp c_2 t^{p/2}, \qquad \|V_t(\varphi)\| \asymp c_3 t^{p/2},$$

*as $t$ goes toward infinity, where $\|\cdot\|$ denotes the Euclidean norm, $p$ is the dimension of $\Theta$ and $V_t^r(\varphi)$, $V_t(\varphi)$ refer here to the sequential importance resampling case, that is, $k_t(\theta, \cdot) = \delta_\theta$.*

The conditions mentioned above amount to assuming that $\pi_t$ is the posterior density of a model regular enough to ensure the existence and asymptotic normality of the maximum likelihood estimator. Under such conditions, $\pi_t$ can be approximated at first order as a Gaussian distribution centered at $\theta_0$ with variance $I(\theta_0)^{-1}/t$, where $I(\theta_0)$ is the Fisher information matrix evaluated at $\theta_0$. The results above are then derived through the Laplace approximation of integrals; see the Appendix. At first glance, it seems paradoxical that $V_t^{\mathrm{sis}}(\varphi)$ converges to zero when $p = 1$. Note, however, that the ratio $V_t(\varphi)/\mathrm{Var}_{\pi_t}(\varphi)$, which measures the precision of the algorithm relative to the variation of the considered function, is likely to diverge even when $p = 1$, since typically $\mathrm{Var}_{\pi_t}(\varphi) \asymp I(\theta_0)^{-1}/t$ as $t \to +\infty$.

That the sequential importance resampling algorithm diverges more quickly than the sequential importance sampling algorithm in this context is unsurprising: when particles are not mutated, the only effect of a selection step is to deplete the particle system. In this respect, we have for any nonconstant function $\varphi$,

$$V_t^{\mathrm{sis}}(\varphi) < V_t^r(\varphi) \le V_t(\varphi).$$

The proof of this inequality is straightforward.

Due to its facility of implementation and the results above, it may be recommended to use the sequential importance sampling algorithm for studying short series of observations, provided that the dimension of $\Theta$ is low. But, in general, one should rather implement a more elaborate particle filter which includes mutation steps in order to counter the particle depletion. A further implication of these results is the following. Consider a dynamical model which involves a fixed parameter $\theta$, and assume that the marginal posterior distributions $\pi(\theta|y_{1:t})$, obtained by marginalizing out latent variables $x_{1:t}$, satisfy the regularity conditions of Theorem 4. Then, following the argument developed in Section 2.5, we get that the rate of divergence of the sequential importance resampling algorithm for this kind of model is at least of order $O(t^{p/2})$, where $p$ is the dimension of this fixed parameter.

3.3. *Sequential importance sampling and resampling for Bayesian filtering and smoothing.* For simplicity we assume that $\pi_t(x_{1:t}) = \pi(x_{1:t}|y_{1:t})$ is the posterior density of a state space model with latent Markov process $(x_t)$, $x_t \in \mathcal{X}$, and observed process $(y_t)$, $y_t \in \mathcal{Y}$, which satisfies the equations

$$y_t|x_t \sim f(y_t|x_t)\, dy_t,$$

$$x_t|x_{t-1} \sim g(x_t|x_{t-1})\, dx_t.$$

We distinguish two types of functions: those which are defined on common dimensions of the spaces $\Theta_t = \mathcal{X}^t$, say, $\varphi : x_{1:t} \to \varphi(x_k)$, for $t \geq k$, and those which are evaluated on the "last" dimension of $\Theta_t$, that is, $\varphi : x_{1:t} \to \varphi(x_t)$. Evaluating these two types of functions amounts to, respectively, "smoothing" or "filtering" the states.

The sequential importance sampling algorithm is usually very inefficient in such a context, whether for smoothing or filtering the states. We illustrate this phenomenon by a simple example. Assume the $t$th mutation step consists of drawing $x_t$ from the prior conditional density $g(x_t|x_{t-1})$, which is usually easy to implement. Consider two evolving particles $\theta_t^{(j)} = x_{1:t}^{(j)}$ with weights $w_t^{(j)}$, $j = 1, 2$. We have

$$\log \frac{w_t^{(1)}}{w_t^{(2)}} = \sum_{k=1}^{t} \log \frac{f(y_k|x_k^{(1)})}{f(y_k|x_k^{(2)})}.$$

Assuming that the joint process $(y_t, x_t^{(1)}, x_t^{(2)})$ is stationary, the sum above typically satisfies some central limit theorem of the form

$$(13) \qquad t^{-1/2} \sum_{k=1}^{t} \log \frac{f(y_k|x_k^{(1)})}{f(y_k|x_k^{(2)})} \xrightarrow{\mathcal{D}} \mathcal{N}(0, \sigma^2),$$

where the limiting distribution is centered for symmetry reasons. Note that this convergence is with respect to the joint probability space of the simulated processes $x_t^{(j)}$, $j = 1, 2$ and the observation process $(y_t)$, while all our previous results were for a given sequence of observations. In this way, (13) yields that the ratio of weights of the two particles either converges or diverges exponentially fast. More generally, when $H$ particles are generated initially, very few of them will have a prominent weight after some iterations, thus leading to very unreliable estimates, whether for smoothing or filtering the states. The algorithm suffers from the curse of dimensionality, in that its degeneracy grows exponentially with the dimension of the space of interest $\Theta_t$.

We now turn to the sequential importance resampling algorithm, and remark first that, for $\varphi : x_{1:t} \to \varphi(x_1)$ and $t > 0$,

$$V_t(\varphi) \geq V_t^r(\varphi) > V_t^{\mathrm{sis}}(\varphi),$$

provided $\varphi$ is not constant. The proof of this inequality is straightforward. The sequential importance resampling algorithm is even more inefficient than the sequential importance sampling algorithm in smoothing the first state $x_1$, because the successive selection steps only worsen the deterioration of the particle system in the $x_1$ dimension. This is consistent with our claim in Section 2.4 that a selection step always degrades the inference on past and current states, but may possibly improve the inference on future states. In this respect, the algorithm is expected to show more capability in filtering the states, and we now turn to the study of the filtering stability.

The functional operator $\mathcal{E}_t$ which appears in the expression for $V_t(\varphi)$, see (9), summarizes two antagonistic effects: on one hand, the weight distortion due to the correction step, and, on the other hand, the rejuvenation of particles due to the application of the kernel $k_t$. Stability will be achieved provided that these two effects compensate in some way.

For simplicity, we assume that the state space $\mathcal{X}$ is included in the real line and that the studied filtering function $\varphi : x_{1:t} \to \varphi(x_t)$ is real-valued. Recall that for the sequential importance resampling algorithm, $k_t$ is given by

$$k_t(x_{1:t-1}^*, dx_{1:t}) = \delta_{x_{1:t-1}^*} q_t(x_t | x_{1:t-1}^*) \, dx_t,$$

for some given conditional probability density $q_t(\cdot|\cdot)$. We assume that $q_t$ only depends on the previous state $x_{t-1}$, and, therefore, defines a Markov transition. The ability of $q_t$ to "forget the past" is usually expressed through its contraction coefficient [see Dobrushin (1956)]

$$\rho_t = \tfrac{1}{2} \sup_{x', x'' \in \mathcal{X}} \|q_t(\cdot | x') - q_t(\cdot | x'')\|_1,$$

where $\|\cdot\|_1$ stands for the $L_1$-norm. Note $\rho_t \leq 1$, and if $\rho_t < 1$, $q_t$ is said to be strictly contractive. Define the variation of a given function $\varphi$ by

$$\Delta \varphi = \sup_{x, x' \in \mathcal{X}} |\varphi(x) - \varphi(x')|.$$

Then the coefficient $\rho_t$ measures the extent to which the application $q_t$ "contracts" the variation of the considered function, that is, for any $x', x'' \in \mathcal{X}$,

$$(14) \qquad \left| \int q_t(x|x')\varphi(x)\,dx - \int q_t(x|x'')\varphi(x)\,dx \right| \le \rho_t \Delta\varphi.$$

Furthermore, it is known [Dobrushin (1956)] that if $q_t$ is such that, for all $x, x', x'' \in \mathcal{X}$,

$$\frac{q_t(x|x')}{q_t(x|x'')} \le C,$$

then its contraction coefficient satisfies $\rho_t \le 1 - C^{-1}$. We therefore make such assumptions in order to prove the stability of the sequential importance resampling algorithm.

THEOREM 5. *Assume that $\Delta\varphi < +\infty$ and there exist constants $C$, $\underline{f}$ and $\bar{f}$ such that, for any $t \ge 0$, $x, x', x'' \in \mathcal{X}$, $y \in \mathcal{Y}$,*

$$(15) \qquad \frac{g(x|x')}{g(x|x'')} \le C, \qquad \frac{q_t(x|x')}{q_t(x|x'')} \le C, \qquad 0 < \underline{f} \le f(y|x) \le \bar{f}.$$

*Then $V_t(\varphi)$ is bounded from above in $t$ (in the sequential importance resampling case).*

This theorem is akin to previous results in the literature [see Del Moral and Guionnet (2001), Le Gland and Oudjane (2004) and most especially, Künsch (2001, 2003)], except that these authors rather consider the stability of some distance (such as the total variation norm of the difference) between the "true" filtering density $\pi_t(x_t)$ and the empirical density computed from the particle system. In fact, Del Moral and Miclo [(2000), page 36] proved that the *actual* variance of the Monte Carlo error is bounded from above over time under similar conditions. Unfortunately, all these results, including ours, require strong assumptions, such as (15), that are unrealistic when $\mathcal{X}$ is not compact. Further research will hopefully provide weaker assumptions, but this may prove an especially arduous problem.

3.4. *Resample-move algorithms, variance estimation.* Following Gilks and Berzuini (2001), we term "resample-move algorithm" any particle filter algorithm which includes an MCMC step in order to reduce degeneracy, as described in Section 2.2. It seems difficult to make general statements about such algorithms and we will rather make informal comments.

The fixed parameter case is especially well behaved. Basic particle filters diverge only at a polynomial rate, as seen in Section 3.2, in contrast with the exponential rate for state-space models. Adding (well-calibrated) MCMC mutation steps should, consequently, lead to stable algorithms in many cases of interest. In fact, it is doubtful that a mutation step must be performed at each iteration to achieve stability. Chopin (2002) argues and provides some experimental evidence that it may be sufficient to perform move steps at a logarithmic rate, that is, the $n$th move step should occur at iteration $t_n \sim \exp(\alpha n)$.

Situations where a latent process intervenes seem less promising. Smoothing the states is especially a difficult problem, and we do not think that there is any solution for circumventing the curse of dimensionality that we have pointed out in the previous section. Even if mutation steps are performed at every iteration, the MCMC transition kernels should themselves suffer from the curse of dimensionality, in that their ability to rejuvenate particles of dimension $t$ is likely to decrease with $t$.

Resample-move algorithms remain an interesting alternative when the considered dynamical model includes a fixed parameter $\theta$. MCMC mutation steps should avoid depletion in simulated values of $\theta$, and make it possible at least to filter the states and estimate the parameter under reasonable periods of time. Unfortunately, the corresponding MCMC transition kernels will often depend on the whole past trajectory, so that long term stability remains uncertain.

In such complicated setups it is necessary to monitor at least numerically the degeneracy of the considered particle filter algorithm. We propose the following method. Run $k$, say $k = 10$, parallel independent particle filters of size $H$. For any quantity to be estimated, compute the average of the $k$ corresponding estimates. This new estimator is clearly consistent and asymptotically normal. Moreover, the computational cost of this strategy is identical to that of a single particle filter of size $kH$, while the obtained precision will be also of the same order of magnitude in both cases, that is to say $\{V_t(\varphi)/(kH)\}^{1/2}$. This method does not, therefore, incur an unnecessary computational load, and allows for assessing the stability of the algorithm through the evolution of the empirical variance of these $k$ estimates.

## APPENDIX

**A.1. Proofs of Theorems 1 and 2.** We start by outlining some basic properties of the sets $\Phi_t^{(d)}$ with respect to linearity. The set $\Phi_t^{(d)}$ is stable through linear transformations, that is, $\varphi \in \Phi_t^{(d)} \Rightarrow M\varphi \in \Phi_t^{(d')}$ if $M$ is a $d' \times d$ matrix of real numbers. In particular, if the vector function $\varphi = (\varphi_1, \ldots, \varphi_d)'$ belongs to $\Phi_t^{(d)}$, then each of its coordinates belongs to $\Phi_t^{(1)}$. The converse proposition is also true. Finally, we have $V_t(M\varphi + \lambda) = MV_t(\varphi)M'$ for any constant $\lambda \in \mathbb{R}^d$, and this relation also holds for the operators $\widetilde{V}_t$ and $\widecheck{V}_t$. Proving these statements is not difficult and is left to the reader.

The proof works by induction with Lemmas A.1–A.3 for Theorem 1, and Lemmas A.1, A.2 and A.4 for Theorem 2. The inductive hypothesis is the following. For a given $t > 0$, it is assumed that for all $\varphi \in \Phi^{(d)}_{t-1}$,

$$(16) \qquad H^{1/2} \left\{ \frac{1}{H} \sum_{j=1}^{H} \varphi(\hat{\theta}^{(j,H)}_{t-1}) - \mathbb{E}_{\pi_{t-1}}(\varphi) \right\} \xrightarrow{\mathcal{D}} \mathcal{N}\{0, \widehat{V}_{t-1}(\varphi)\}.$$

LEMMA A.1 (Mutation). *Under the inductive hypothesis, we have*

$$H^{1/2} \left\{ \frac{1}{H} \sum_{j=1}^{H} \psi(\theta^{(j,H)}_t) - \mathbb{E}_{\tilde{\pi}_t}(\psi) \right\} \xrightarrow{\mathcal{D}} \mathcal{N}\{0, \widetilde{V}_t(\psi)\}$$

*for any measurable* $\psi \colon \Theta_t \to \mathbb{R}^d$ *such that the function* $\mu \colon \theta_{t-1} \mapsto \mathbb{E}_{k_t(\theta_{t-1}, \cdot)}\{\psi(\cdot) - \mathbb{E}_{\tilde{\pi}_t}(\psi)\}$ *belongs to* $\Phi^{(d)}_{t-1}$ *and there exists* $\delta > 0$ *such that* $\mathbb{E}_{\tilde{\pi}_t} \|\psi\|^{2+\delta} < +\infty$.

PROOF. We assume that $\psi$ is real-valued $(d = 1)$. The generalization to $d > 1$ follows directly from the Cramér–Wold theorem and the linearity properties stated above.

Let $\bar{\psi} = \psi - \mathbb{E}_{\tilde{\pi}_t}(\psi)$, $\mu(\theta_{t-1}) = \mathbb{E}_{k_t(\theta_{t-1}, \cdot)}\{\bar{\psi}(\cdot)\}$, $\sigma^2(\theta_{t-1}) = \text{Var}_{k_t(\theta_{t-1}, \cdot)}\{\bar{\psi}(\cdot)\}$ and $\sigma_0^2 = \mathbb{E}_{\pi_{t-1}}(\sigma^2)$. We have $\mathbb{E}_{\pi_{t-1}}(\mu) = 0$, and by Jensen's inequality,

$$\sigma_0^2 = \mathbb{E}_{\pi_{t-1}}[\text{Var}_{k_t(\theta_{t-1}, \cdot)}\{\psi(\cdot)\}] \leq \mathbb{E}_{\pi_{t-1}}[\mathbb{E}_{k_t(\theta_{t-1}, \cdot)}\{\psi(\cdot)^2\}]$$
$$\leq \{\mathbb{E}_{\tilde{\pi}_t} |\psi|^{(2+\delta)}\}^{2/(2+\delta)} < +\infty,$$

which makes it possible to apply the law of large numbers for particle filters to $\sigma^2$,

$$(17) \qquad H^{-1} \sum_{j=1}^{H} \sigma^2(\theta^{(j,H)}_{t-1}) \to \sigma_0^2 \qquad \text{almost surely.}$$

Defining

$$(18) \qquad \nu(\theta_{t-1}) = \mathbb{E}_{k_t(\theta_{t-1}, \cdot)}\{|\bar{\psi}(\cdot) - \mu(\theta_{t-1})|^{2+\delta}\}$$

$$(19) \qquad \leq 2^{1+\delta}\{\mathbb{E}_{k_{t-1}(\theta_{t-1}, \cdot)}|\bar{\psi}(\cdot)|^{2+\delta} + |\mathbb{E}_{k_{t-1}(\theta_{t-1}, \cdot)}\bar{\psi}(\cdot)|^{2+\delta}\}$$

$$(20) \qquad \leq 2^{2+\delta}\{\mathbb{E}_{k_{t-1}(\theta_{t-1}, \cdot)}|\bar{\psi}(\cdot)|^{2+\delta}\},$$

where (19) comes from the $C_r$ inequality and (20) from Jensen's inequality, we deduce that

$$\mathbb{E}_{\pi_{t-1}}(\nu) \leq 2^{2+\delta}\mathbb{E}_{\tilde{\pi}_t}|\psi|^{2+\delta} < +\infty.$$

This inequality ensures that the expectations defining $\nu$ in (18) (and, similarly, those defining $\mu$ and $\sigma^2$) are finite for almost every $\theta_{t-1}$. It follows that

$$H^{-1} \sum_{j=1}^{H} \nu(\theta_{t-1}^{(j,H)}) \to \mathbb{E}_{\pi_{t-1}}(\nu) \qquad \text{almost surely,}$$

and combining this result with (17), we obtain the almost sure convergence of

$$\begin{aligned}
(21) \quad \rho_H &= \frac{\sum_{j=1}^{H} \nu(\theta_{t-1}^{(j,H)})}{\{\sum_{j=1}^{H} \sigma^2(\theta_{t-1}^{(j,H)})\}^{(2+\delta)/2}} \\
&= H^{-\delta/2} \frac{H^{-1} \sum_{j=1}^{H} \nu(\theta_{t-1}^{(j,H)})}{\{H^{-1} \sum_{j=1}^{H} \sigma^2(\theta_{t-1}^{(j,H)})\}^{(2+\delta)/2}} \to 0.
\end{aligned}$$

Let $T_H = H^{-1/2} \sum_{j=1}^{H} \bar{\psi}(\theta_t^{(j,H)})$, $S_{t-1}$ denote the sigma-field generated by the random variables forming the triangular array $(\hat{\theta}_{t-1}^{(j,H)})_{j \leq H}$, that is, the particle system at time $t-1$, and $\mu_H = \mathbb{E}(T_H | S_{t-1})$. Conditional on $S_{t-1}$, the $\bar{\psi}(\theta_t^{(j,H)})$'s form a triangular array of independent variables which satisfy the Liapunov condition, see (21), and have variances whose mean converges to $\sigma_0^2$, see (17). Therefore [Billingsley (1995), page 362], the following central limit theorem for triangular arrays of independent variables holds:

$$(22) \qquad (T_H - \mu_H) | S_{t-1} \xrightarrow{\mathcal{D}} \mathcal{N}(0, \sigma_0^2).$$

Since $\mathbb{E}_{\pi_{t-1}}(\mu) = 0$ and $\mu \in \Phi_{t-1}^{(d)}$, we have also, by applying (16) to the function $\mu$,

$$(23) \qquad \mu_H = H^{-1/2} \sum_{j=1}^{H} \mu(\hat{\theta}_{t-1}^{(j,H)}) \xrightarrow{\mathcal{D}} \mathcal{N}\{0, \widehat{V}_{t-1}(\mu)\}.$$

The characteristic function of $T_H$ is

$$\begin{aligned}
\Phi_{T_H}(u) &= \mathbb{E}\{\exp(iuT_H)\} \\
&= \mathbb{E}[\exp(iu\mu_H)\mathbb{E}\{\exp(iuT_H - iu\mu_H)|S_{t-1}\}],
\end{aligned}$$

where $\mathbb{E}\{\exp(iuT_H - iu\mu_H)|S_{t-1}\}$ is the characteristic function of $T_H - \mu_H$ conditional on $S_{t-1}$, which according to (22) converges to $\exp(-\sigma_0^2 u^2/2)$. It follows from (23) that

$$\exp(iu\mu_H)\mathbb{E}\{\exp(iuT_H - iu\mu_H)|S_{t-1}\} \xrightarrow{\mathcal{D}} \exp(-\sigma_0^2 u^2/2 + iuZ),$$

where $Z$ is a random variable distributed according to $\mathcal{N}\{0, \widehat{V}_{t-1}(\mu)\}$. The expectation of the left-hand side term converges to the expectation of the right-hand side term following the dominated convergence theorem, and this completes the proof. $\square$

LEMMA A.2 (Correction). *Let $\varphi \in \Phi_t^{(d)}$, assume the inductive hypothesis holds and the function $\theta_t \mapsto 1$ belongs to $\Phi_t^{(1)}$. Then*

$$H^{1/2}\left\{\frac{\sum_{j=1}^H w_t^{(j,H)}\varphi(\theta_t^{(j,H)})}{\sum_{j=1}^H w_t^{(j,H)}} - \mathbb{E}_{\pi_t}(\varphi)\right\} \xrightarrow{\mathcal{D}} \mathcal{N}\{0, V_t(\varphi)\}.$$

PROOF. Let $\bar{\varphi} = \varphi - \mathbb{E}_{\pi_t}(\varphi)$. For notational convenience we assume that $d = 1$, but the generalization to $d \geq 1$ is straightforward. It is clear that the vector function $\psi = (v_t \cdot \bar{\varphi}, v_t)'$ fulfills the conditions mentioned in Lemma A.1, and as such satisfies

$$H^{1/2}\left\{\frac{1}{H}\sum_{j=1}^H \begin{pmatrix} v_t(\theta_t^{(j,H)})\bar{\varphi}(\theta_t^{(j,H)}) \\ v_t(\theta_t^{(j,H)}) \end{pmatrix} - \begin{pmatrix} 0_{\mathbb{R}^d} \\ 1 \end{pmatrix}\right\} \xrightarrow{\mathcal{D}} \mathcal{N}\{0, \widetilde{V}_t(\psi)\}.$$

Then, resorting to the $\delta$-method with function $g(x, y) = x/y$, we obtain

$$H^{1/2}\frac{\sum_{j=1}^H v_t(\theta_t^{(j,H)})\bar{\varphi}(\theta_t^{(j,H)})}{\sum_{j=1}^H v_t(\theta_t^{(j,H)})} \xrightarrow{\mathcal{D}} \mathcal{N}(0, \mathcal{V}),$$

where $\mathcal{V} = \{(\partial g/\partial x, \partial g/\partial y)(0, 1)\}\widetilde{V}_t(\psi)\{(\partial g/\partial x, \partial g/\partial y)(0, 1)\}' = \widetilde{V}_t\{v_t \cdot (\varphi - \mathbb{E}_{\pi_t}\varphi)\}$. The left-hand side term is unchanged if we replace the $v_t(\theta_t^{(j,H)})$'s by the weights $w_t^{(j,H)}$, since they are proportional. $\square$

LEMMA A.3 (Selection, multinomial resampling). *Let $\widehat{V}_t(\varphi) = V_t(\varphi) + \mathrm{Var}_{\pi_t}(\varphi)$ and assume the particle system is resampled according to the multinomial scheme. Then, under the same conditions as in Lemma* A.2,

$$H^{1/2}\left\{\frac{1}{H}\sum_{j=1}^H \varphi(\hat{\theta}_t^{(j,H)}) - \mathbb{E}_{\pi_t}(\varphi)\right\} \xrightarrow{\mathcal{D}} \mathcal{N}\{0, \widehat{V}_t(\varphi)\}.$$

PROOF. The proof is similar to that of Lemma A.1. Assume $d = 1$, denote by $S_t$ the sigma-field generated by the random variables $(\theta_t^{(j,H)}, w_t^{(j,H)})_{j \leq H}$ and let $\bar{\varphi} = \varphi - \mathbb{E}_{\pi_t}(\varphi)$, $T_H = H^{-1/2}\sum_{j=1}^H \bar{\varphi}(\hat{\theta}_t^{(j,H)})$ and $\mu_H = \mathbb{E}(T_H|S_t)$. Conditional on $S_t$, $T_H$ is, up to a factor $H^{-1/2}$, a sum of independent draws from the multinomial distribution which produces $\bar{\varphi}(\theta_t^{(j,H)})$ with probability $w_t^{(j,H)}/\sum_{j=1}^H w_t^{(j,H)}$. Then, as in Lemma A.1, we have

$$(T_H - \mu_H)|S_t \xrightarrow{\mathcal{D}} \mathcal{N}(0, \sigma_0^2),$$

where this time $\sigma_0^2 = \mathrm{Var}_{\pi_t}(\varphi)$, which is the limit as $H \to +\infty$ of the variance of the multinomial distribution mentioned above. The proof is completed along the same lines as in Lemma A.1.    $\square$

LEMMA A.4 (Selection, residual resampling).    *Let $\widehat{V}_t(\varphi)$ take the value given by* (7) *and assume the particle system is resampled according to the residual resampling scheme. Then, under the same conditions as in Lemma* A.2,

$$H^{1/2}\left\{\frac{1}{H}\sum_{j=1}^{H}\varphi(\hat{\theta}_t^{(j,H)}) - \mathbb{E}_{\pi_t}(\varphi)\right\} \xrightarrow{\mathcal{D}} \mathcal{N}\{0, \widehat{V}_t(\varphi)\}.$$

PROOF.    The proof is identical to that of Lemma A.2, except that conditional on $S_t$, $T_H$ is $H^{-1/2}$ times a constant, plus a sum of independent draws from the multinomial distribution described in Section 2.1. This yields a different value for $\sigma_0^2$,

$$\sigma_0^2 = \mathbb{E}_{\tilde{\pi}_t}\{r(v_t) \cdot \varphi^2\} - \frac{1}{\mathbb{E}_{\tilde{\pi}_t}\{r(v_t)\}}[\mathbb{E}_{\tilde{\pi}_t}\{r(v_t) \cdot \varphi\}]^2.$$

In addition, we also have to make sure that the number of these independent draws $H^r$ tends toward infinity. In fact, $H^r/H \to \mathbb{E}_{\tilde{\pi}_t}[r(\nu_t)]$. To see this, consider

$$H^r/H - H^{-1}\sum_{j=1}^{H}r\{v_t(\theta_t^{(j,H)})\} = H^{-1}\sum_{j=1}^{H}[r(H\rho_j) - r\{v_t(\theta_t^{(j,H)})\}],$$

where $H\rho_j = v_t(\theta_t^{(j,H)})/\{H^{-1}\sum_j v_t(\theta_t^{(j,H)})\}$, see Section 2.1, so that the difference above should eventually be zero as $H^{-1}\sum_j v_t(\theta_t^{(j,H)}) \to 1$. More precisely, we have $|r(x) - r(y)| \le 1$, in general, and $r(x) - r(y) = x - y$ provided $|x - y| < \varepsilon$ and $r(x) \in [\varepsilon, 1 - \varepsilon]$ for any $\varepsilon < 1/2$. Therefore, assuming that $\{H^{-1}\sum_j v_t(\theta_t^{(j,H)})\}^{-1} \in [1 - \varepsilon', 1 + \varepsilon']$ for some $\varepsilon' > 0$ and $H$ large enough, we get that the sum above should be zero plus something bounded from above by the proportion of particles such that $\varepsilon' v_t(\cdot) > 1/2$ or $r\{v_t(\cdot)\} \notin [\varepsilon' v_t(\cdot), 1 - \varepsilon' v_t(\cdot)]$. This proportion can be made as small as necessary.    $\square$

**A.2. Proof of Theorem 3.**    Let $\varphi : \Xi_t \to \mathbb{R}^d$ and $\bar{\varphi} = \varphi - \mathbb{E}_{\pi_t}(\varphi) = \varphi - \mathbb{E}_{\pi_t^m}(\varphi)$ for a given $t \ge 0$. To simplify notation, it is assumed that $d = 1$, but the adaptation to the general case is straightforward. All quantities related to the "marginalized" particle filter are distinguished by the $m$-suffix. For instance, $\mathcal{E}_t^m(\varphi)$ stands for the function $\xi_t \mapsto \mathbb{E}_{k_t^m(\xi_t, \cdot)}\{v_t^m(\cdot)\varphi(\cdot)\}$, in agreement with the definition of $\mathcal{E}_t(\varphi)$ in (10). In this respect, the marginal weight function $v_t^m(\cdot)$ is $\tilde{\pi}_t^m(\cdot)/\pi_t^m(\cdot)$, and if we define the "conditional" weight function $v_t^c(\lambda_t|\xi_t) = \pi_t^c(\lambda_t|\xi_t)/\tilde{\pi}_t^c(\lambda_t|\xi_t)$, we have the identity

$$v_t(\theta_t) = v_t^m(\xi_t) v_t^c(\lambda_t|\xi_t).$$

It follows from (12) that

$$\mathbb{E}_{\pi_{t-1}^c}\{\mathcal{E}_t(\bar{\varphi})\} = \mathbb{E}_{k^m}\{v_t^m \bar{\varphi} \mathbb{E}_{\tilde{\pi}_t^c}(v_t^c)\} = \mathcal{E}_t^m(\bar{\varphi}),$$

since $\mathbb{E}_{\tilde{\pi}_t^c}(v_t^c) = 1$, and by induction, we show similarly, for $k \le t$, that

$$\mathbb{E}_{\pi_k^c}\{\mathcal{E}_{k+1:t}(\bar{\varphi})\} = \mathcal{E}_{k+1:t}^m(\bar{\varphi}).$$

Hence, for $k \le t$,

$$\mathbb{E}_{\tilde{\pi}_k}[\{v_k \mathcal{E}_{k+1:t}(\bar{\varphi})\}^2] = \mathbb{E}_{\tilde{\pi}_k^m}[(v_k^m)^2 \mathbb{E}_{\tilde{\pi}_k^c}\{v_k^c \mathcal{E}_{k+1:t}\bar{\varphi}\}^2]$$
$$\ge \mathbb{E}_{\tilde{\pi}_k^m}[\{v_k^m \cdot \mathcal{E}_{k+1:t}^m(\bar{\varphi})\}^2],$$

by Jensen's inequality. From the closed form (9) of $V_t(\varphi)$, we deduce the inequality $V_t^m(\varphi) \le V_t(\varphi)$ for the case when the selection step follows the multinomial scheme. Alternatively, if the selection step consists of residual resampling, let $\underline{\varphi} = \varphi - \mathbb{E}_{\tilde{\pi}_t}\{r(v_t)\varphi\}/\mathbb{E}_{\tilde{\pi}_t}\{r(v_t)\}$. Then

$$R_t(\varphi) - R_t^m(\varphi) = \mathbb{E}_{\tilde{\pi}_t}\{r(v_t)\underline{\varphi}^2\} - \mathbb{E}_{\tilde{\pi}_t^m}\{r(v_t^m)\underline{\varphi}^2\} + \frac{\{\mathbb{E}_{\tilde{\pi}_t^m}r(v_t^m)\underline{\varphi}\}^2}{\mathbb{E}_{\tilde{\pi}_t^m}r(v_t^m)}$$
$$\ge \mathbb{E}_{\tilde{\pi}_t^m}[\{\mathbb{E}_{\tilde{\pi}_t^c}r(v_t) - r(v_t^m)\}\underline{\varphi}^2],$$

and since $\mathbb{E}_{\tilde{\pi}_t^c}(v_t) = v_t^m$, we have $\mathbb{E}_{\tilde{\pi}_t^c}\lfloor v_t \rfloor \le \lfloor v_t^m \rfloor$, hence $\mathbb{E}_{\tilde{\pi}_t^c}r(v_t) \ge r(v_t^m)$, and, consequently, $R_t(\varphi) \ge R_t^m(\varphi)$ for any $\varphi$. It is then easy to show by induction that the desired inequality is also verified in the residual case.

**A.3. Regularity conditions and proof of Theorem 4.** Let $\pi_0(\theta)$ denote the prior density and $p(y_{1:t}|\theta)$ the likelihood of the $t$ first observations, so that through Bayes formula,

$$\pi_t(\theta) = \pi(\theta|y_{1:t}) \propto \pi_0(\theta) p(y_{1:t}|\theta).$$

Let $l_t(\theta) = \log p(y_{1:t}|\theta)$. The following statements are assumed to hold almost surely:

1. The maximum $\hat{\theta}_t$ of $l_t(\theta)$ exists and converges as $t \to +\infty$ to $\theta_0$ such that $\pi_0(\theta_0) > 0$ and $\tilde{\pi}_0(\theta_0) > 0$.
2. The matrix

$$\Sigma_t = -\left\{\frac{1}{t}\frac{\partial^2 l_t(\theta)}{\partial\theta\partial\theta'}\right\}^{-1}$$

is positive definite and converges to $I(\theta_0)$, the Fisher information matrix at $\theta_0$.

3. There exists $\Delta > 0$ such that

$$0 < \delta < \Delta \implies \limsup_{t \to +\infty} \left[ \frac{1}{t} \sup_{\|\theta - \hat{\theta}_t\| > \delta} \{ l_t(\theta) - l_t(\hat{\theta}_t) \} \right] < 0.$$

4. The functions $\pi_0(\theta)$ and $l_t(\theta)$ are six-times continuously differentiable, the partial derivatives of order six of $l_t(\theta)/t$ are bounded on any compact set $\Theta' \subset \Theta$, and the bound does not depend on $t$ and the observations.
5. $\varphi : \Theta \to \mathbb{R}^d$ is six-times continuously differentiable, $\varphi'(\theta_0) \neq 0$.

For convenience, we start with the one-dimensional case ($p = 1$). The Laplace approximation of an integral [see, e.g., Tierney, Kass and Kadane (1989)] is

$$\int \psi(\theta) \exp\{-th(\theta)\} \, d\theta$$

$$= (2\pi/t)^{1/2} \sigma \exp\{-t\hat{h}\}$$

$$\times [\hat{\psi} + \tfrac{1}{2}\{\sigma^2 \hat{\psi}'' - \sigma^4 \hat{\psi}' \hat{h}''' + \tfrac{5}{12}\sigma^6 \hat{\psi} \hat{h}''' - \tfrac{1}{4}\sigma^4 \hat{\psi} \hat{h}^{iv}\} t^{-1} + O(t^{-2})],$$

where hats on $\psi$, $h$ and their derivatives indicate evaluation at the point which minimizes $h$, and $\sigma = -(1/\hat{h}'')^{1/2}$. This approximation remains valid for a function $h_t$ depending on $t$, provided that the fluctuations of $h_t$ or its derivatives can be controlled in some way. Conditions above allow, for instance, for applying this approximation to the functions $h_{t,1} = -l_t(\theta)/t$ and $h_{t,2}(\theta) = -2l_t(\theta)/t$; see Schervish [(1995), page 446] for technical details. It is necessary, however, to assume that $\psi(\theta_0) \neq 0$, so that $\psi$ is either strictly positive or strictly negative at least in a neighborhood of $\theta_0$. Since $V_t^{\mathrm{sis}}(\varphi) = V_t^{\mathrm{sis}}(\varphi + \lambda)$ for any $\lambda \in \mathbb{R}$, we assume without loss of generality that $\varphi(\theta_0) \neq 0$. $V_t^{\mathrm{sis}}(\varphi)$ equals

$$(24) \quad \begin{aligned} &\frac{\int \psi_1(\theta) p(y_{1:t}|\theta)^2 \, d\theta - 2\mathbb{E}_{\pi_t}(\varphi) \int \psi_2(\theta) p(y_{1:t}|\theta)^2 \, d\theta}{\{\int \pi(\theta) p(y_{1:t}|\theta) \, d\theta\}^2} \\ &+ \frac{\{\mathbb{E}_{\pi_t}(\varphi)\}^2 \int \psi_3(\theta) p(y_{1:t}|\theta)^2 \, d\theta}{\{\int \pi(\theta) p(y_{1:t}|\theta) \, d\theta\}^2}, \end{aligned}$$

where $\psi_1 = \pi_0(\theta)^2 \varphi(\theta)^2 / \tilde{\pi}_0(\theta)$, $\psi_2 = \pi_0(\theta)^2 \varphi(\theta) / \tilde{\pi}_0(\theta)$ and $\psi_3 = \pi_0(\theta)^2 / \tilde{\pi}_0(\theta)$. Combining the appropriate Laplace approximations, we get that

$$V_t^{\mathrm{sis}}(\varphi) = \frac{t^{1/2}}{2(\pi \Sigma_t)^{1/2}}$$

$$\times \frac{[\psi_1(\hat{\theta}_t) - 2\mathbb{E}_{\pi_t}(\varphi)\psi_2(\hat{\theta}_t) + \{\mathbb{E}_{\pi_t}(\varphi)\}^2 \psi_3(\hat{\theta}_t) + At^{-1} + O(t^{-2})]}{\{\pi_0(\hat{\theta}_t) + Bt^{-1} + O(t^{-2})\}^2}$$

$$= \frac{t^{1/2}}{2(\pi \Sigma_t)^{1/2}} \frac{\{\varphi(\hat{\theta}_t) - \mathbb{E}_{\pi_t}(\varphi)\}^2 + A\tilde{\pi}_0(\hat{\theta}_t)\pi_0(\hat{\theta}_t)^{-2} t^{-1} + O(t^{-2})}{\tilde{\pi}_0(\hat{\theta}_t)\{1 + B\pi_0(\hat{\theta}_t)^{-1} t^{-1} + O(t^{-2})\}^2},$$

where $A$ is the sum of $O(t^{-1})$ terms corresponding to the three Laplace expansions of the numerator, and $B$ is the $O(t^{-1})$ term of the denominator. Since $\varphi(\hat{\theta}_t) - \mathbb{E}_{\pi_t}(\varphi) = O(t^{-1})$, $\Sigma_t = I(\theta_0) + O(t^{-1})$ and $\psi(\hat{\theta}_t) = \psi(\theta_0) + O(t^{-1})$ for any continuous function $\psi$, we get through appropriate derivations that

$$V_t^{\text{sis}}(\varphi) = \frac{I(\theta_0)^{1/2}\varphi'(\theta_0)^2}{2\pi^{1/2}\tilde{\pi}_0(\theta_0)} t^{-1/2} + O(t^{-3/2}).$$

Derivations in multidimensional cases are much the same, except that notation is more cumbersome. When $p > 1$, the factor $t^{-1/2}$ in the Laplace expansion is replaced by $t^{-p/2}$, so that in the ratio (24) we get a factor $t^{p/2}$, and since the $t^{p/2}$ terms cancel as in the one-dimensional case, the actual rate of divergence is $t^{p/2-1}$, and this completes the first part of the proof.

In the sequential importance resampling case (multinomial scheme), $q_t(\theta, \cdot) = \delta_\theta$ and $\tilde{\pi}_t = \pi_{t-1}$, and according to (9),

$$(25) \qquad V_t(\varphi) = V_t^{\text{sis}}(\varphi) + \sum_{k=1}^{t} \mathbb{E}_{\pi_{k-1}}\left[\frac{\pi_t}{\pi_{k-1}}\{\varphi - \mathbb{E}_{\pi_t}(\varphi)\}\right]^2.$$

Then through a direct adaptation of expansions above we obtain a divergence rate for $V_t(\varphi)$ of order $(\sum_{k=0}^{t}(t-k)^{p/2-1}) = O(t^{p/2})$. For the residual case, it follows from (11) and (25) that

$$V_t^r(\varphi) = V_t^{\text{sis}}(\varphi) + \sum_{k=1}^{t} R_{k-1}\left[\frac{\pi_t}{\pi_{k-1}}\{\varphi - \mathbb{E}_{\pi_t}(\varphi)\}\right].$$

The difficulty in this case is that the noncontinuous function $r(\cdot)$ takes part in the expression for $R_k(\cdot)$, see (8). It is clear, however, that the Laplace expansion can be generalized to cases where regularity conditions for the likelihood and other functions are fulfilled only locally around $\theta_0$. The additional assumption that $\pi_t(\theta_0)/\pi_{t-1}(\theta_0)$ is not an integer for any $t > 0$ allows $r(v_t)$ to be six-times continuously differentiable in a neighborhood around $\theta_0$, and, therefore, makes it possible to expand the terms of the sum above, which leads to a rate of divergence of order $O(t^{p/2})$ in the same way as in the multinomial case.

**A.4. Proof of Theorem 5.** As a preliminary, we state without proof the following inequality. Let $\varphi, \psi : \mathbb{R} \to \mathbb{R}$ such that $\varphi \geq 0$, $\sup \psi \geq 0$ and $\inf \psi \leq 0$. Then

$$(26) \qquad \Delta(\varphi\psi) \leq \sup \varphi \cdot \Delta\psi.$$

Due to particular cancelations, the weight function $v_t(x_{1:t})$ only depends on $x_{t-1}$ and $x_t$ in the state space case

$$(27) \qquad v_t(x_{1:t}) = v_t(x_{t-1}, x_t) \propto \frac{f(y_t|x_t)g(x_t|x_{t-1})}{q_t(x_t|x_{t-1})}.$$

Straightforward consequences of this expression are the identities

$$(28) \qquad \pi_t(x_t|x_{t-1}) = \frac{q_t(x_t|x_{t-1})v_t(x_{t-1},x_t)}{\int q_t(x|x_{t-1})v_t(x_{t-1},x)\,dx},$$

$$(29) \qquad \pi_{t+1}(x_{t+1}|x_k) = \frac{\int \pi_t(x_t|x_k)q_{t+1}(x_{t+1}|x_t)v_{t+1}(x_t,x_{t+1})\,dx_t}{\int \pi_t(x_t|x_k)q_{t+1}(x|x_t)v_{t+1}(x_t,x)\,dx_t\,dx},$$

for $k < t$, where $\pi_t(x_t|x_k)$ denotes the conditional posterior density of $x_t$ given $x_k$ and the $t$ first observations, that is, $\pi_t(x_t|x_k) = \pi(x_t|x_k, y_{1:t}) = \pi(x_t|x_k, y_{k+1:t})$. We start by proving some useful lemmas.

LEMMA A.5. *The conditional posterior density $\pi_t(x_t|x_k)$, $k < t$, defines a Markov transition from $x_k$ to $x_t$ whose contraction coefficient is less than or equal to $(1 - C^{-2})^{t-k}$.*

PROOF. This is adapted from Künsch (2001). For $x_k, x'_k, x_{k+1} \in \mathcal{X}$, $k < t$,

$$\frac{\pi_t(x_{k+1}|x_k)}{\pi_t(x_{k+1}|x'_k)} = \frac{g(x_{k+1}|x_k)p(y_{k+1:t}|x'_k)}{g(x_{k+1}|x'_k)p(y_{k+1:t}|x_k)} \leq C^2,$$

since $g(x_{k+1}|x_k) \leq Cg(x_{k+1}|x'_k)$ and

$$p(y_{k+1:t}|x'_k) = \int g(x_{k+1}|x'_k)p(y_{k+1:t}|x_{k+1})\,dx_{k+1}$$

$$\leq C \int g(x_{k+1}|x_k)p(y_{k+1:t}|x_{k+1})\,dx_{k+1}.$$

Therefore, the contraction coefficients of Markov transitions $\pi_t(x_{k+1}|x_k)$ and $\pi_t(x_t|x_k)$ are less than or equal to, respectively, $(1 - C^{-2})$ and $(1 - C^{-2})^{t-k}$. □

LEMMA A.6. *Let $\lambda$ be a probability density on $\mathcal{X}$ and $h(x|x')$ a conditional probability density defining a Markov transition on $\mathcal{X}$. Then for any $x' \in \mathcal{X}$, $y \in \mathcal{Y}$,*

$$\frac{\int f(y|x)h(x|x')\,dx}{\mathbb{E}_{\lambda(x'')}\{\int f(y|x)h(x|x'')\,dx\}} \leq 1 + \rho_h C_f,$$

*where $\rho_h$ is the contraction coefficient of $h(\cdot|\cdot)$, and $C_f = \bar{f}/\underline{f} - 1$.*

PROOF. It follows from the definition of $\rho_h$ [see (14)] that for $x', x'' \in \mathcal{X}$,

$$\left| \int f(y|x)h(x|x')\,dx - \int f(y|x)h(x|x'')\,dx \right| \leq \rho_h(\bar{f} - \underline{f})$$

and therefore,

$$\sup_{x' \in \mathcal{X}} \left\{ \int f(y|x)h(x|x')\,dx \right\} \leq \mathbb{E}_{\lambda(x'')} \left\{ \int f(y|x)h(x|x'')\,dx \right\} + \rho_h(\bar{f} - \underline{f}),$$

so that

$$\frac{\sup_{x' \in \mathcal{X}}\{\int f(y|x)h(x|x')\,dx\}}{\mathbb{E}_{\lambda(x'')}\{\int f(y|x)h(x|x'')\,dx\}} \leq 1 + \rho_h \frac{(\bar{f} - \underline{f})}{\mathbb{E}_{\lambda(x'')}\{\int f(y|x)h(x|x'')\,dx\}}$$

$$\leq 1 + \rho_h \left(\frac{\bar{f}}{\underline{f}} - 1\right). \qquad \square$$

LEMMA A.7. *Let* $\rho = 1 - C^{-1}$ *and* $\rho_2 = 1 - C^{-2}$. *Then for* $k < t$,

$$\Delta\mathcal{E}_{k+1:t}\{\varphi - \mathbb{E}_{\pi_t}(\varphi)\} \leq \prod_{i=1}^{t-k}(1 + \rho\rho_2^{i-1}C_f)\rho_2^{t-k}\Delta\varphi,$$

*for any real-valued filtering function,* $\varphi: x_{1:t} \to \varphi(x_t)$.

PROOF. Let $\bar{\varphi} = \varphi - \mathbb{E}_{\pi_t}(\varphi)$. Note the arguments of $\mathcal{E}_{k+1:t}(\bar{\varphi})$ are $x_{1:k}$ in general, but in the case considered in Section 3.3 it only depends on $x_k$ and is therefore treated as a function $\mathcal{X} \to \mathcal{X}$. For the sake of clarity, we treat the case $k = t - 2$, but the reasoning is easily generalized. The following decomposition is deduced from identity (28):

$$\mathcal{E}_{t-1:t}(\bar{\varphi})(x_{t-2})$$
$$= \mathbb{E}_{q_{t-1}(x_{t-1}|x_{t-2})}\{\upsilon_{t-1}(x_{t-2}, x_{t-1})\mathcal{E}_t(\bar{\varphi})(x_{t-1})\}$$
$$= \mathbb{E}_{q_{t-1}(x_{t-1}|x_{t-2})}\{\upsilon_{t-1}(x_{t-2}, x_{t-1})\}\mathbb{E}_{\pi_{t-1}(x_{t-1}|x_{t-2})}\{\mathcal{E}_t(\bar{\varphi})(x_{t-1})\}.$$

It follows from (27) that the first term satisfies

$$\mathbb{E}_{q_{t-1}(x_{t-1}|x_{t-2})}\{\upsilon_{t-1}(x_{t-2}, x_{t-1})\} \propto \int f(y_{t-1}|x_{t-1})g(x_{t-1}|x_{t-2})\,dx_{t-1},$$

where the proportionality constant can be retrieved by remarking that the expectation of this term with respect to $\pi_{t-2}$ equals one and, therefore,

$$\mathbb{E}_{q_{t-1}(x_{t-1}|x_{t-2})}\{\upsilon_{t-1}(x_{t-2}, x_{t-1})\}$$
$$= \frac{\int f(y_{t-1}|x_{t-1})g(x_{t-1}|x_{t-2})\,dx_{t-1}}{\mathbb{E}_{\pi_{t-2}(x_{t-2})}\{\int f(y_{t-1}|x_{t-1})g(x_{t-1}|x_{t-2})\,dx_{t-1}\}}$$
$$\leq 1 + \rho C_f$$

according to Lemma A.6. Note $\pi_{t-2}(x_{t-2})$ denotes the $\pi_{t-2}$-marginal density of $x_{t-2}$. It follows from the decomposition above and the inequality in (26) that

$$\Delta\mathcal{E}_{t-1:t}(\bar{\varphi}) \leq (1 + \rho C_f)\Delta\psi,$$

where $\psi$ is the function

$$\psi(x_{t-2}) = \mathbb{E}_{\pi_{t-1}(x_{t-1}|x_{t-2})}\{\mathcal{E}_t(\bar{\varphi})(x_{t-1})\}$$
$$= \mathbb{E}_{\pi_{t-1}(x_{t-1}|x_{t-2})}[\mathbb{E}_{q_t(x_t|x_{t-1})}\{\upsilon_t(x_{t-1}, x_t)\bar{\varphi}(x_t)\}].$$

Note that $\psi$ does take positive and negative values, since the expectation of $\mathcal{E}_{t-1\,:\,t}(\bar\varphi)$ with respect to $\pi_{t-2}$ is null. We now decompose $\psi$ in the same way,

$$\psi(x_{t-2}) = \mathbb{E}_{\pi_{t-1}(x_{t-1}|x_{t-2})}[\mathbb{E}_{q_t(x_t|x_{t-1})}\{\upsilon_t(x_{t-1}, x_t)\}]\mathbb{E}_{\pi_t(x_t|x_{t-2})}\{\bar\varphi(x_t)\},$$

by consequence of the identity (29). The expectation of the first term with respect to $\pi_{t-1}(x_{t-2})$ equals one, so that

$$\mathbb{E}_{\pi_{t-1}(x_{t-1}|x_{t-2})}[\mathbb{E}_{q_t(x_t|x_{t-1})}\{\upsilon_t(x_{t-1}, x_t)\}]$$

$$= \frac{\int \pi_{t-1}(x_{t-1}|x_{t-2})f(y_t|x_t)g(x_t|x_{t-1})\,dx_{t-1}\,dx_t}{\mathbb{E}_{\pi_{t-1}(x_{t-2})}\{\int \pi_{t-1}(x_{t-1}|x_{t-2})f(y_t|x_t)g(x_t|x_{t-1})\,dx_{t-1}\,dx_t\}}$$

$$\leq 1 + \rho\rho_2 C_f,$$

according to Lemmas A.5 and A.6. Resorting again to inequality (26), we get

$$\Delta\psi \leq (1 + \rho\rho_2 C_f)\rho_2^2 \Delta\varphi,$$

which leads to the desired inequality, and this completes the proof of Lemma A.7. □

To conclude the proof of Theorem 5, remark that $\mathbb{E}_{\bar\pi_k}(\upsilon_k) = 1$. Therefore,

$$\upsilon_k(x_{k-1}, x_k) = \frac{f(y_k|x_k)g(x_k|x_{k-1})/q_k(x_k|x_{k-1})}{\mathbb{E}_{\bar\pi_k(x_{1\,:\,k})}\{f(y_k|x_k)g(x_k|x_{k-1})/q_k(x_k|x_{k-1})\}}$$

$$\leq C^2 \bar{f}/\underline{f},$$

and since the expectation of the function $\mathcal{E}_{k+1\,:\,t}\{\varphi - \mathbb{E}_{\pi_t}(\varphi)\}$ with respect to $\pi_k$ is null, the function $\mathcal{E}_{k+1\,:\,t}\{\varphi - \mathbb{E}_{\pi_t}(\varphi)\}$ is ensured to take positive and negative values, so that

$$\sup_{x_k\in\mathcal{X}} |\mathcal{E}_{k+1\,:\,t}\{\varphi - \mathbb{E}_{\pi_t}(\varphi)\}(x_k)| \leq \Delta\mathcal{E}_{k+1\,:\,t}\{\varphi - \mathbb{E}_{\pi_t}(\varphi)\}$$

and, finally,

$$\mathbb{E}_{\bar\pi_k}[\upsilon_k^2 \mathcal{E}_{k+1\,:\,t}\{\varphi - \mathbb{E}_{\pi_t}(\varphi)\}^2]$$

$$\leq C^4(\bar{f}/\underline{f})^2 \prod_{i=1}^{t-k}(1 + \rho\rho_2^{i-1}C_f)^2 \rho_2^{2(t-k)}(\Delta\varphi)^2$$

$$\leq C^4(\bar{f}/\underline{f})^2 \exp\left(2\rho C_f \sum_{i=1}^{t-k}\rho_2^{i-1}\right)\rho_2^{2(t-k)}(\Delta\varphi)^2$$

$$\leq C^4(\bar{f}/\underline{f})^2 \exp\{2\rho C_f/(1-\rho_2)\}\rho_2^{2(t-k)}(\Delta\varphi)^2.$$

It follows from (9) that $V_t(\varphi)$ is bounded from above by a convergent series.

**Acknowledgments.** This paper is the fourth part of my Ph.D. thesis, defended on March 2003 at Université Pierre et Marie Curie, Paris. It was inspired in part by Hans Künsch's lectures on particle filters at the "Summer School on Advanced Computational Methods for Statistical Inference," Luminy, September 2001, and benefited from helpful comments by Christian Robert, Hans Künsch, Eric Moulines, Pierre Del Moral, one anonymous referee and an Associate Editor.

DEPARTMENT OF MATHEMATICS
UNIVERSITY OF BRISTOL
UNIVERSITY WALK
BRISTOL BS8 1TW
UNITED KINGDOM
E-MAIL: nicolas.chopin@bristol.ac.uk